\input{BoxedEPS.tex}
\SetTexturesEPSFSpecial
\HideDisplacementBoxes
\documentclass[11pt]{amsproc}
\usepackage{amsmath,amsfonts,amssymb}

\usepackage{graphicx}

\newcommand{\Z}{\mathbb Z}

\newcommand{\G}{\Gamma}

\newtheorem{theorem}{Theorem}
\newtheorem{proposition}[theorem]{Proposition}
\newtheorem{lemma}[theorem]{Lemma}
\newtheorem{corollary}[theorem]{Corollary}

\newtheorem{question}{Question}

\def\endproof{~\hfill${\blacksquare}$\medskip}
\def\sd{\rtimes}
\def\parno{\par\noindent}
\begin{document}

\title{Metabelian  Wreath Products are LERF}

\author{Roger C. Alperin}
\address{Department of Mathematics, San Jose State University, San Jose, CA 95192}
\email{alperin@math.sjsu.edu  }
\maketitle 
\section{Introduction}
The subgroup $S$ of $\G$ is separable means that $S$ is closed in the profinite topology of $\G$.  A  finitely generated (f.g.) group is LERF if all its
f.g. subgroups are separable.

Gruenberg's
theorem \cite{G} asserts that the wreath product
$A\wr Q$ is residually finite (r.f.) iff $A$ and $Q$ are r.f. and either $A$ is abelian or $Q$ is finite.
We seek general conditions which will characterize the wreath products which are LERF. Since LERF entails
r.f. Gruenberg's theorem gives a first restriction. Since subgroups of LERF groups are LERF we shall
 assume $A$ and $Q$ are LERF. 

In the  case where $Q$ is finite and $A$ is a finitely generated LERF group then 
up to finite index $\G=A\wr Q$ is direct product  of finitely many LERF groups. 
\begin{question} For which f.g. groups  $X$ is  $X^n$ LERF for all $n\ge 1$.
\end{question}
\parno  Polycyclic groups are LERF by a  theorem of Malcev and thus  $X$  polycyclic gives a
positve answer to the question.
By a theorem of M. Hall the free group $F_n$ of rank $n$ is LERF while Mihailova has shown that $F_2\times F_2$
is not
\cite{L-S}; so any group
$X$ which contains
$F_2$ can  not give a positive answer to this question. 

Our main result concerns the other case when $A$ is f.g. abelian; it generalizes the recent result of 
\cite{deC}: $A\wr \Z$ is LERF for any f.g. abelian group. We show that when $A, Q$ are f.g. abelian
then $\G=A\wr Q$ is LERF. Using the Magnus embedding it follows that free metabelian groups are LERF.
\section{Lemma}

 Here is small modification of 
an important lemma from \cite{deC}. We say a subgroup $S$ of  $\G$ is `strongly
separable'  if all  the finite index subgroups of $S$ are separable  in $\G$.

\begin{lemma}\label{de} Let $\G$
be f.g and r.f. Suppose $f: \G\rightarrow Q$ is a surjective  homomorphism with abelian kernel $K$.    If  $S$ is a
f.g. subgroup such that $f(S)$ is of finite index in $Q$ then $S$ is strongly separable in $\G$. 
\end{lemma}
\proof   If
$S'$ is a subgroup of finite index in
$S$ then its image $Q'$ is  finite index in $Q$.  We can show $S'$ is separable in $\G$ using
$f':f^{-1}(Q')\rightarrow Q'$ with kernel $K$.  If $S'$ is closed in a subgroup of finite index in $\G$
 then it is closed in
$\G$ so
$S$ is strongly closed. Also, if $f(S)=Q'$ of finite index in $Q$ then $f^{-1}(Q')$ is a 
subgroup of finite index in $\G$; since separable in a subgroup, 
say $f^{-1}(Q')$, of finite index in $\G$ is separable in $\G$, it suffices to prove the case when $f(S)=Q$. 
We may suppose then
that $S'=S$ and $f(S)=Q$.
It now follows that  $\G=KS$; consequently  $S \cap K$ is normal in $\G$.  
If we show
$S/S\cap K$ is separable in
$\G/S\cap K$ then since $S$ is f.g we can also separate  $S$ in $\G$. Thus we need only consider the group $\bar{\G}=\G/S\cap K$
which is a split extension; we  assume then $S\cap K=\{1\}$ and  consequently $\G=K\sd S$. 

It now suffices to show that the profinite closure $\bar{S}$ of
$S$ in $\G$ meets  $K$ trivially, $\bar{S}\cap K=\{1\}$.  For then, if $x\in \bar{S}$ then $f(x)x^{-1}\in \bar{S}\cap K=\{1\}$;
hence
$x\in S$ so $S$ is closed. 
Consider then $x\in \bar{S}\cap K$; if $x\ne 1$ then  since $\G$ is r.f. we can
choose a normal subgroup $N$ of finite index in $\G$ so that $x\notin N$; 
consider $L=N\cap K$ and $T=N\cap S$; then $x\notin L$. Since $T$ normalizes $L$ then $M=\bigcap_{t\in S\ mod\
T}L^s$ is finite index in $L$, and $x\notin M$; since $M\sd T$ is finite index in $\G$ and does not contain $x$, then $x\notin
\bar{S}$; thus we must have  $x=1$.
\endproof

\section{Main Results}
We now assume $A$,  $Q$ are f.g. abelian.  
Any subgroup $R$ of $Q$ is contained as a subgroup of finite index in $R_1$ and there is a (retract)
homomorphism
$\pi: Q\rightarrow R_1$  so that
$\pi(r)=r$, $r\in R$. We can in fact achieve the case that $\pi$ is split and  that $Q=R_1\times R_2$ with
$R_1,R_2$ infinite, when $R$ is non-trivial, not of finite index and  $Q$ has rank greater than 1. Then 

$$\G=A\wr Q =A^{R_1\times R_2}\ltimes (R_1\times R_2) $$
$$=(A^{R_2\times R_1}\ltimes R_1)\times (A^{R_1\times R_2}\ltimes R_2) $$
$$=(A^{R_2}\wr R_1)\times (A^{R_1}\wr R_2) $$
$$=\G^{(2)}_{R_1}\times \G^{(1)}_{R_2}.$$

\begin{proposition} Supose that $A, Q$  are f.g. abelian  then $\G^{(2)}_{R_1}, \G^{(1)}_{R_2}$ are f.g and  r.f. 
\end{proposition}
\proof Since $\G$ is f.g.  these  quotient groups are also f.g. By Gruenberg's theorem these are also r.f. since
any direct sum of cyclics is r.f.
\endproof

\begin{theorem} Suppose  $A$ and $Q$ are f.g.  abelian
groups then  $A\wr Q$ is LERF.
\end{theorem}
\proof  

Let
$K=A^Q$,
$\G=A\wr Q$; we have the natural homomorphism
$f:
\G\rightarrow Q$.    
Let $x\notin S$, a f.g. subgroup of $\G$; $R=f(S)$. We may assume that $R$ is not of finite index in $Q$ for
then the Theorem follows from the  Lemma. Also the case of $Q$ having rank one is proven in \cite{deC}. As in the
remarks above we consider after passing to a subgroup of finite index in $S$ if necessary,  a retract $R_1$ with
$R_1\times R_2$ a subgroup of finite index in $Q$, each
$R_i$ infinite. Since a subgroup $S$ is separable if a subgroup of finite index is separable in a subgroup of
finite index in $\G$, we may now assume that $R_1\times R_2=Q$.

If
$x$ has a non-trivial image in $R_2$  we project $\G$ to $\G^{(1)}_{R_2}$. In this group the image of $S$ is
trivial and the image of
$x$ is not so we can map to a finite quotient to separate because  
$\G^{(1)}_{R_2}$ is r.f. by the Proposition.

Suppose then $x$ has trivial image in $R_2$, then we can replace $\G$ by $\G^{(2)}_{R_1}$ and consider there 
$x\notin S$; we can separate them
since $S$ is closed in 
$\G^{(2)}_{R_1}$ by the Lemma since $f(S)$ is of finite index in  $R_1$ and  since   
$\G^{(2)}_{R_1}$ is f.g and  r.f.  by the Proposition. 
\endproof

\begin{corollary} Free metabelian groups are LERF.
\end{corollary}
\proof Via the Magnus embedding, \cite{M}, the free metabelian group can be embedded in a group $A\wr Q$ where $A$
and
$Q$ are f.g.  abelian groups. The corollary follows now since subgroups of LERF groups are LERF. 
\endproof

\end{document}